\documentclass[a4paper,11pt]{amsart}

\usepackage{amsmath,amssymb,amsfonts,latexsym}

\newtheorem{theorem}{Theorem}[section]
\newtheorem{lemma}[theorem]{Lemma}
\newtheorem{corollary}[theorem]{Corollary}

\begin{document}

\title[Attractor Dimension of Rotating NSE-$\alpha$]{Attractor Dimensions of
  Three-Dimensional
  Navier-Stokes-$\alpha$ Model for Fast Rotating Fluids on Generic-Period Domains: 
  Comparison with Navier-Stokes Equations}

\author{Bong-Sik Kim}
\address{Department of Mathematics \& Natural Sciences,
 American University of Ras Al Khaimah,  
United Arab Emirates}
\email{bkim@aurak.ac.ae}

%%%%%%%%%%%%%%%%%%%%%%%%%%%%%%%%%%%%%%%%%%%%%%%%%%%%%%%%%%%%%%%%%%%

\begin{abstract}
  The three-dimensional Navier-Stokes-$\alpha$ model for fast rotating 
   geophysical fluids 
   is considered. The Navier-Stokes-$\alpha$ model is a nonlinear dispersive regularization of 
   the exact Navier-Stokes equations obtained by Lagrangian averaging and 
    tend to 
    the Navier-Stokes equations as $\alpha\rightarrow 0^+$.
    We estimate upper bounds for the dimensions of global attractors 
    and 
     study  the dependence of the dimensions on the parameter $\alpha$.  
    All the estimates are uniform in $\alpha$, and our estimate of attractor dimensions
      remain finite when $\alpha\rightarrow 0^+$.
\vspace{2mm}


\vspace{2mm}


\end{abstract}

%\thanks{This work was supported by ..}

%%%%%%%%%%%%%%%%%%%%%%%%%%%%%%%%%%%%%%%%%%%%%%%%%%%%%%%%%%%%%%%%%%%

\maketitle

%%%%%%%%%%%%%%%%%%%%%%%%%%%%%%%%%%%%%%%%%%%%%%%%%%%%%%%%%%%%%%%%%%%

\section{Introduction}

%%%%%%%%%%%%%%%%%%%%%%%%%%%%%%%%%%%%%%%%%%%%%%%%%%%%%%%%%%%%%%%%%%%

We consider the three-dimensional rotating Navier-Stokes-$\alpha$ equations (RNS-$\alpha$)
with periodic boundary conditions in a torus $\mathbb{T}^3=[0, 2\pi a_1]\times[0, 2\pi a_2]\times [0, 2\pi a_3]$:

\begin{eqnarray}
  \frac{\partial v}{\partial t}+(u\cdot\nabla)v+v_j\nabla u^j
         +\Omega e_3\times u   & = & -\nabla p+\nu\Delta v+f \label{eq1}\\ 
             \nabla\cdot v = \nabla\cdot u  =  0\  \ & \mbox{and}& v(t,x)|_{t=0}=v_0 \nonumber  \\
             u &=& (I-\alpha^2\Delta )^{-1}v  \nonumber
\end{eqnarray}
where $x=(x_1,x_2,x_3)\in \mathbb{T}^3$, $v = v(t,x) = (v_1(t,x), v_2(t,x), v_3(t,x)$ 
is the velocity field, $p=p(x,t)$ is the pressure of a homogeneous incompressible fluid, 
$\nu$ is the viscosity, and $f=f(x)$ is a divergence free body force.
$\Omega$ is the Coriolis parameter, which is twice the angular velocity of the rotation around 
the vertical unit vector $e_3=(0,0,1)$.
 The system  (\ref{eq1}) reduces to the exact rotating Navier-Stokes equations (RNS)
  when $\alpha \rightarrow  0^+$. 

Kim and Nicolaenko \cite{Kim2} established the existence and global regularity of solutions of the system
  (\ref{eq1}) and proved the existence of its global attractor.  
In this paper, we estimate the dimension of a global attractor of the system (\ref{eq1})
and give special attention to the limiting case when $\alpha\rightarrow 0^+$, that is,
when RNS-$\alpha$ equations tend to the RNS equations. 
We focus on generic-period domains and eliminate 
 \textit{nontrivial resonant parts} (strict three-wave resonant interactions), 
 which are essentially related  to Rossby wave in physics.
 General periodic cases have to deal with
 nontrivial resonant parts, which will be covered in a separate article.

  Ilyin and Titi \cite{Ilyin1} estimated attractor dimensions for two-dimensional navier-Stokes-$\alpha$
  equations. Their estimates, however, blow up as $\alpha\rightarrow 0^+$. 
  Gibbon and Holm \cite{Gibbon1} obtained length-scale estimates for NS-$\alpha$ equations in 
  terms of the Reynolds number, which blow up in the limit when $\alpha\rightarrow 0^+$, too. 
  Several other time-averaged estimates related 
  to NS-$\alpha$ equations don't remain finite in the Navier-Stokes limit,
  except the cases where there are known equivalent upper bounds
  for the Navier-Stokes equations (see  Table 1 in \cite{Gibbon1}). 
  They analyzed the system in the context of the filtered velocity $u = (I-\alpha^2\Delta )v$.
  Instead,  we study the system from the perspective of 
   the non-filtered velocity $v = (I-\alpha^2\Delta )^{-1}u$.
  The Helmholtz inverse operator $\mathcal{R}_\alpha = (I-\alpha^2\Delta )^{-1}$ plays a  crucial role
   in the process, leading to 
  uniform estimates that remain finite in the Navier-Stokes limit as $\alpha\rightarrow 0^+$.
  
  We consider domain parameters, $a_1, a_2$ and $a_3$
 bounded away from both zero and infinity.
$v$ and $u$ are vector fields defined on
$D=\mathbb{T}^3\times [0, T]$
for any $T>0$. Periodicity of the boundary conditions leads naturally to a Fourier
 representation of the fields, that is
\[
  v = \sum_{n} v_n e^{i (n_1x_1/a_1+n_2x_2/a_2+n_3/a_3)} =
    \sum_{n} v_n e^{i\check{n}\cdot x},
\]
where $v_n$'s are Fourier Coefficients,
$n=(n_1,n_2,n_3)\in \mathbb{Z}^3$ , and
 $\check{n}=(\check{n}_1,\check{n}_2,\check{n}_3)$ 
are  wave numbers with $
 \check{n}_j=n_j/a_j$ for $j=1,2,3$.
We set $a_1=1$ without loss of generality and
define the Fourier-Sobolev space of
divergence free periodic vector fields as follows:
{\small
\[  
   H^{s}=\left\{ v\in [L^2(\mathbb{T}^3)]^3 \ |\
           v=\sum_{n\in \mathbb{Z}^3}v_ne^{i\check{n}\cdot x},\,
                      v^{*}_n=v_{-n},\ v_0=0,\ \check{n}\cdot v_n=0,\
                      ||v||_{s}^2<\infty\right\},    
\]
}
with the norm
            \[ ||v||_{s}^2=
             \sum_{n\in \mathbb{Z}^3}|\check{n}|^{2s}|v_n|^2.
            \]
Here $v_n^*$ is the complex conjugate of $v_n$.
The corresponding inner product is denoted by $<\cdot ,\cdot >_s$. 
We set 
$H^0=H$ when $s=0$. Also, $<\cdot ,\cdot >_0=<\cdot ,\cdot >$, $||\cdot ||_0=|\cdot |$, and
$||\cdot ||_1=||\cdot ||$.
We assume that
\[
  \int_{\mathbb{T}^3}v(x,0)\,dx = 0 \ \ \ \mbox{ and }\ \ \
  \int_{\mathbb{T}^3}f(x)\,dx = 0\ \ \ \mbox{ for all } t\geq 0.
\]
This yields  $\int_{\mathbb{T}^3}v(x,t)\,dx = 0$ for all $t\geq 0$,  and
allows for the use of the Poincar\'{e} inequality.

We denote $P_L$ as the usual Leray projection onto the divergence free subspace and 
introduce the Helmholtz inverse operator, $\mathcal{R}_\alpha = (I-\alpha^2\Delta )^{-1}$,
 which is given by
\[ \mathcal{R}_\alpha v = (I-\alpha^2\triangle )^{-1}v.
\]
A bilinear operator $B_\alpha$ on divergence free vector fields \cite{Kim2} is define by
\[ B_\alpha(u,v) = P_L\left[ (\mathcal{R}_\alpha u\cdot\nabla )v+v_j\nabla (\mathcal{R}_\alpha u)_j\right] 
   = -P_L\left[\mathcal{R}_\alpha v\times \mathrm{curl}\,v\right].
\]
This bilinear operator has a connection with the classical Navier-Stokes bilinear operator
\[ B(u,v)=P_L[(u\cdot\nabla )v].\]
\begin{lemma}\label{bilinear-relation}
   For every $u,v,w\in H^1$
\begin{equation*}
  < B_{\alpha}(u,v),\ w> =
   <B(\mathcal{R}_{\alpha}u,v),\ w>-
   <B(w,v),\ \mathcal{R}_{\alpha}u>
\end{equation*}
\end{lemma}
\noindent\textbf{Proof.} 
  \begin{eqnarray*}
    <B_{\alpha}(u,v),w>&=& <-P_L[\mathcal{R}_{\alpha}u\times (\nabla\times v)], \ w> \\
                                   &=& <P_L[(\mathcal{R}_{\alpha}u\cdot\nabla )v+\sum_{j=1}^3v_j\cdot\nabla (\mathcal{R}_{\alpha}u)_j],\ w>\\
                                   &=& <P_L[(\mathcal{R}_{\alpha}u\cdot\nabla )v],\ w>+<P_L[\sum_{j=1}^3v_j\cdot\nabla (\mathcal{R}_{\alpha}u)_j],\ w> \\
                                   &=& <B(\mathcal{R}_{\alpha}u,v),\ w>+<\sum_{j=1}^3v_j\cdot\nabla (\mathcal{R}_{\alpha}u)_j,\ w>\\
                                   &=& <B(\mathcal{R}_{\alpha}u,v),\ w>-<B(w,v),\ \mathcal{R}_{\alpha}u>
  \end{eqnarray*}
 For the second equality we use the identity, $(a\cdot\nabla )b
 =\nabla(a\cdot b)-(b\cdot\nabla )a
 -a\times\mathrm{curl}\,b
 -b\times\mathrm{curl}\,a$, to get 
\begin{equation*}
   \nabla (v\cdot\mathcal{R}_{\alpha}u)-\mathcal{R}_{\alpha}u\times\mbox{curl}\,v=
   (\mathcal{R}_{\alpha}u\cdot\nabla )v+(v\cdot\nabla )\mathcal{R}_{\alpha}u+v\times\mathrm{curl}\,\mathcal{R}_{\alpha}u.
\end{equation*}
Noticing that $(v\cdot\nabla )\mathcal{R}_{\alpha}u+v\times\mathrm{curl}\,\mathcal{R}_{\alpha}u=
  \sum_{j=1}^{3}v_j\nabla (\mathcal{R}_{\alpha}u)_j$, we can get the equality. For the last equality, we directly calculate the second inner product
 such as
\begin{eqnarray*}
  <\sum_{j=1}^{3}v_j\nabla (\mathcal{R}_{\alpha}u)_j,\ w>&=&\sum_{i,j=1}^{3}\int_Qv_j\frac{\partial (\mathcal{R}_{\alpha}u)_j}{\partial x_i}w_i\,dx \\
   &=& \sum_{i,j}\int_Qw_i \frac{\partial (\mathcal{R}_{\alpha}u)_j}{\partial x_i}v_j\,dx \\
   &=& -\sum_{i,j}\int_Qw_i\frac{\partial v_j}{\partial x_i}(\mathcal{R}_{\alpha}u)_j\,dx \\
   &=& -<(w\cdot\nabla )v,\ \mathcal{R}_{\alpha}u>\\
   &=&-<P_L[(w\cdot\nabla )v],\ \mathcal{R}_{\alpha}u>\\
   &=&<-B(w,v),\ \mathcal{R}_{\alpha}u>.
\end{eqnarray*}
Therefore, the result follows.\ \ \ \ $\blacksquare$
\medskip 

Now we rewrite
Eq(\ref{eq1}) in terms of the unfiltered velocity $v$: 
\begin{equation}\label{eq2}
    \frac{\partial v}{\partial t}+
   \Omega\,P_LJP_L\mathcal{R}_{\alpha}v + \nu Av+B_{\alpha}(v,v)
    = f,
\end{equation}
where   $A=-P_L\Delta$ is the Stokes operator and $J$
 is a rotation matrix given by
 \[ J = \left( \begin{array}{ccc}
                  0 & -1 & 0 \\
          1 & 0 & 0 \\
          0 & 0 & 0
        \end{array} \right) .
 \]
In  Fourier-Sobolev space, the action
$P_L$ on $n$-th Fourier component of a vector field is given by
  $P_Lv=\sum_{n}(P_nv_n)e^{i\check{n}\cdot x}$ and
  $P_nv_n=
  (v_n-\frac{\check{n}\cdot v_n}{|\check{n}|^2}
   \check{n})$ with
   \[ P_n=I-\frac{1}{|\check{n}|^2}
        \left( \begin{array}{ccc}
               n_1^2 & \frac{n_1n_2}{a_2} & \frac{n_1n_3}{a_3} \\
              \frac{n_1n_2}{a_2} & \frac{n_2^2}{a_2^2}& \frac{n_2n_3}{a_2a_3} \\
              \frac{n_1n_3}{a_3} & \frac{n_2n_3}{a_2a_3} & \frac{n_3^2}{a_3^2}
               \end{array}
        \right).
   \]
 The Helmholtz inverse operator $\mathcal{R}_{\alpha}$ commutes with curl and, for 
 each wave number $n$, 
 \[ (\mathcal{R}_{\alpha})_n = \frac{1}{1+\alpha^2|\check{n}|^2}. \]
Then, for each wave number $n\in\mathbb{Z}^3$, the  RNS-$\alpha$ equations have the form
 \begin{equation*}\label{eq3}
  \frac{\partial v_n}{\partial t}+\frac{1}{1+\alpha^2|\check{n}|^2}
   \Omega\,P_nJP_n v_n+\nu |\check{n}|^2v_n
    +B_{\alpha}(v,v)_n = f_n,
\end{equation*}
where
\begin{equation*}\label{eq4}
  B_{\alpha}(v,v)_n=
   -iP_n\sum_{k+m=n}
   \frac{1}{1+\alpha^2|\check{k}|^2}(v_k\times
   (\check{m}\times v_m)).
\end{equation*}
The existence of unique regular solutions for all $\Omega$
greater than some threshold $\Omega_0$ has been proved
in \cite{BMN2} for $\alpha = 0$ and in  \cite{Kim2} for $\alpha > 0$.
\begin{theorem}[\cite{BMN2}] \label{Th1}
   For every triplet of positive real numbers $(a_1, a_2, a_3)$, the following result holds.
  Let $s>1/2$ and $v_0\in H^s(\mathbb{T}^3)$ a divergence-free vector field.
  Then there exits a constant $\Omega_0>0$, depending on $||v_0||_{s}$, $||f||_{s-1}, \nu,$
  and the domain parameter
   $(a_1, a_2, a_3)$, such that for all $\Omega\geq \Omega_0$,
  there is a unique global solution
  \[ v(t)\in C([0,\infty):H^s(\mathbb{T}^3))
    \cap L^2((0,\infty): H^{s+1}(\mathbb{T}^3)) 
   \]
   to the three-dimensional rotating Navier-Stokes equations ($\alpha = 0$ in (\ref{eq1})).
   Furthermore, if $f$ is independent of $t$, then 
   there exists a global attractor for the three-dimensional rotating
   Navier-Stokes equations bounded in $H^s$; such an attractor has a finite fractal dimension
   and attracts every weak Leray solution as $t\rightarrow +\infty$.
 \end{theorem}
  
 \begin{theorem}[\cite{Kim2}] \label{Th2} 
 For every triplet of positive real numbers $(a_1, a_2, a_3)$, the following result holds.
   Let $s>5/2$ and $v_0\in H^s(\mathbb{T}^3)$ a divergence-free vector field.
  Then there exits a constant $\Omega_{\alpha}>0$, depending on $||v_0||_{s}$,
  $||f||_{s-1}, \nu,$ and the domain parameter
   $(a_1, a_2, a_3)$,
  such that for all $\Omega\geq \Omega_{\alpha}$,
  there is a unique global solution
  \[ v(t)\in C([0,\infty):H^s(\mathbb{T}^3))
    \cap L^2((0,\infty): H^{s+1}(\mathbb{T}^3)) 
   \]
   to the equation for any $\alpha \geq 0$.
   Moreover, all the estimates are uniform in $\alpha$
   (i.e., 
   the estimates don't blow up as $\alpha\rightarrow 0^+$).  
   If $f$ is independent of $t$, then 
   there exists a global attractor for the system (\ref{eq1})
   bounded in $H^s$ and its fractal dimension is finite.
 \end{theorem}
 
{\bf Remark:} The solutions of the three-dimensional rotating Naiver-Stokes-$\alpha$ equations 
 uniformly converge in $L^2$ to those 
 of the three-dimensional rotating Navier-Stokes equations 
 as $\alpha\rightarrow 0^+$ (see Section 8 of \cite{Kim2}).
\smallskip

\noindent We consider the Eq. (\ref{eq1}) in the limit as $\Omega\rightarrow +\infty$,
which gives
 \emph{resonant limit $\alpha$-equations}. Working with the resonant limit $\alpha$-equations 
on specific periodic domains (generic periods),
 we obtain upper estimates of dimensions of global
attractors for 
the resonant limit $\alpha$-equations, which approximate the dimensions of the global attractors
for three-dimensional RNS-$\alpha$ equations on generic periodic domains:

\begin{theorem}\label{main} {\bf (Main Result)}.\ 
  Let $\mathcal{A}_{\alpha}$ be a global attractor of the Eq.(\ref{eq4}) in $H$. 
  Then its Hausdorff dimension $d_H(\mathcal{A}_\alpha)$ and the fractal dimension
  $d_F(\mathcal{A}_\alpha )$  are finite and 
  satisfy, for an absolute constant $K_\alpha$, the estimate
  \[  d_H(\mathcal{A}_\alpha) < K_\alpha\left(\frac{\rho_{V}}{\nu_\alpha}\right)^{2}
  \ \ \mbox{and} \ \ d_F(\mathcal{A}_\alpha )\leq 2d_H(\mathcal{A}_\alpha)
  \]
  where $c_l, \tilde{c}, c_1$, and $d$ are absolute constants, 
  $K(\alpha) =  \left(\frac{c_l}{d}\right)^{3/2} 
         (24c(\alpha)+1)^{3/2} \tilde{c}^{3/2}$,
           $ c^2(\alpha) = \frac{1}{1+\alpha^2c_1}\left[\frac{1}{1+\alpha^2c_1}+\frac{1}{\alpha^2c_1}\right]$,
           $\rho^2_V = 2|f|^2/(\nu^2\lambda_1^2)$, and $\lambda_1$ the first eigenvalue
           of $A=-P_L\Delta$. 
     In particular, for $\alpha =0$, the global attractor has a sharp upper bound
     \[  d_H(\mathcal{A}_0) < \tilde{K}\left(\frac{\rho_{V}}{\nu_0}\right)^{6/5}
     \]
\end{theorem}
\noindent Observe that $\lim_{\alpha\rightarrow 0^+} K(\alpha) =  
            \left(\frac{(24\sqrt{2}+1)c_l\tilde{c}}{d}\right)^{3/2} \equiv K_0 <\infty$
             so that
  the estimates doesn't blow up when $\alpha\rightarrow 0^+$. 
  In general,  $K_0 \neq \tilde{K}$.   
  
  \smallskip
  
Accepting the point of view that the dimension of a global attractor for Navier-Stokes equations
is associated with the number of degree of freedom in turbulent
flows \cite{Constantin1}, then these finite dimension estimates gives a rigorous justification
that asymptotic dynamics of turbulent rotating fluids can be described by
  two-dimensional and three-component (2D-3C) non-steady Navier-Stokes equations 
  when $\Omega$ is large enough.
Also, the attractor dimensions are related to the fundamental length scale $\ell$ of 
turbulent flows, below which wave interactions do not affect its dynamics. For $\alpha$ equations, 
the fundamental length scale in terms of the attractor dimension $d_F(\mathcal{A_\alpha})$ is
      \[ \ell_\alpha^2 \thicksim \left( \frac{|\mathbb{T}^3|}{d_F(\mathcal{A}_{\alpha})}\right)
              \leq \ell_0^2 \thicksim \left( \frac{|\mathbb{T}^3|}{d_F(\mathcal{A}_{0})}\right).
      \]

 As shown in \cite{BMN1, BMN2, Kim2}, the estimates 
 depend crucially on the period
 of the torus, $a_1, a_2, a_3$. If we omit the nonlinearity in (\ref{eq1}) and
  the viscous and forcing terms, we end up with
  the system 
  \begin{equation*}
  \partial_t v + (I-\alpha^2\Delta )^{-1}\Omega\,Jv = -\nabla p,\ \ 
  \nabla \cdot v = 0,
\end{equation*}
 which describes the propagation of waves, called \emph{Poincar\'{e} waves} or \emph{inertial waves}.
 The corresponding dispersion law relating the pulsation $\omega$ to the wavenumber
 $\xi\in\mathbb{R}^3$ is
 \[ \omega (\xi) = \pm \Omega\frac{\xi_3}{|\xi|}. 
 \]
 Therefore, the two-dimensional part of the initial data evolves according to
 two-dimensional Euler or Navier-Stokes equations, and the three-dimensional part
 generates waves, which propagate very rapidly in the domain with a speed of $\Omega$.
 Chemin \emph{et al.} \cite{Chemin1} detailed well about the propagation of high-speed waves in
 the fluid, and we brief their explanation here.  
 The time average of these waves vanishes, but they carry a non-zero energy.
 The wavenumbers of these waves are bounded as $\Omega\rightarrow \infty$ and a priori no short
 wavelengths are created. On periodic flows, Poincar\'{e} waves
 persist for long times, and interact not only with the limit two-dimensional flow,
 but also with themselves.
 A wave $\xi$ interacts with a wave $\xi^\prime$ and generates another wave $\xi^{\prime\prime}$ 
 provided 
 \begin{equation*}\label{chemin}
   \xi + \xi^\prime = \xi^{\prime\prime}\ \ \mbox{and}\ \ 
    \frac{\xi_3}{|\tilde{\xi}|}+\frac{\xi^\prime_3}{|\tilde{\xi}^\prime|}
      = \frac{\xi^{\prime\prime}_3}{|\tilde{\xi}^{\prime\prime}|},\ \ 
      \mbox{where}\ \  \tilde{\xi} = \left( \frac{\xi_1}{a_1},\frac{\xi_2}{a_2},\frac{\xi_3}{a_3}\right) , 
  \end{equation*}
 which are the usual resonance conditions in the three-wave interaction problem.
 In the periodic case, all the components of $\xi$, $\xi^\prime$ and $\xi^{\prime\prime}$ are integers,
 and the above conditions turn out to be Diophantine equations which do not have integer solutions for almost 
 all $(a_1,a_2,a_3)$ except the trivial solutions given by symmetries.
 Therefore, generically in the sizes of the periodic box, the waves
 do not interact with themselves and only interact with the two-dimensional underlying flow.
 Mathematically to handle these resonant wave interactions,
  Babin, Mahalov and Nicolaenko \cite{BMN1, BMN2} introduced
  the Poincar\'{e} group operator (see Section \ref{limit-equations}).
  Then they utilized methods of small denominators and Diophantine 
  incommensurability conditions on the domain geometrical parameters $a_1, a_2, a_3$
  to investigate the fast singular oscillating limits of Eq.(\ref{eq1}) as $\Omega\rightarrow \infty$.
  In that approach, the collective contribution to the dynamics 
  made by fast Poincar\'{e} waves is accounted for by rigorous
  estimates of wave resonances and quasi-resonances via small divisor analysis. 
  We start off  the next section by applying such an approach to derive the resonant limit $\alpha$-equations.

%%%%%%%%%%%%%%%%%%%%%
 \section{Resonant  limit $\alpha$-equations}\label{limit-equations}
%%%%%%%%%%%%%%%%%%%%%%%%
  
 Poincar\'{e} propagator $E_{\alpha}(\Omega t) = e^{\Omega\,t\,P_LJP_L\mathcal{R}_{\alpha}}$
 is defined as the unitary group solution $E_{\alpha}(-\Omega\,t)\Phi_0
  =\Phi (t)\ $($E_{\alpha}(0)=I$ is the identity ) to the linear Poincar\'{e} problem:
\[
  \partial_t \Phi + \Omega\,P_LJP_L\mathcal{R}_{\alpha}\Phi = 0,\ \ \Phi|_{t=0}=\Phi_0,\ \mbox{with}\  
  \nabla \cdot \Phi_0 = 0.
\]
Denote $M_{\alpha}=P_LJP_L\mathcal{R}_{\alpha}$ and
$M_{\alpha n}=(M_{\alpha})_n=
     \frac{1}{1+\alpha^2|\check{n}|^2}P_nJP_n$ for each wavenumber $n$. 
 The matrix $M_{\alpha n}$
has the eigenvalues, $\pm i\omega_{\alpha}(n)$, where
\begin{equation*}\label{trans9}
   \omega_{\alpha}(n) = \frac{\check{n}_3}{(1+\alpha^2|\check{n}|^2)|\check{n}|} = \omega_{\alpha n},
   \ \ |\check{n}| =
    \sqrt{\theta_1 n_1^2+\theta_2 n_2^2+\theta_3 n_3^2},
   \ \ \theta_j = \frac{1}{a_j^2}.
\end{equation*}
In Fourier space,
\begin{eqnarray*}
  E_{\alpha}(\Omega\,t)_n &=& 
          \cos (\Omega\,\omega_{\alpha n}t)I
           +\frac{1}{|\check{n}|}\sin (\Omega\,\omega_{\alpha n}t)R_n,
           \label{trans10}\\
        &=& \frac{1}{2}\left[
        e^{i\Omega\omega_{\alpha n}t}\left( I-i\frac{1}{|\check{n}|}R_n\right)+
        e^{-i\Omega\omega_{\alpha n}t}\left( I+i\frac{1}{|\check{n}|}R_n\right)\right],
        \nonumber
\end{eqnarray*}
where the matrix $iR_n$ is the Fourier transform of the curl
vector; $(\mathrm{curl}\,v)_n=iR_{n}v_n=
i\check{n}\times v_n$ with
   \[ R_n=
        \left( \begin{array}{ccc}
                 0 & -\check{n}_3 & \check{n}_2 \\
              \check{n}_3 & 0 & -n_1 \\
              -\check{n}_2 & n_1 & 0
               \end{array}
        \right) .
   \]
   
  Next, we set $V(t):= E_\alpha (\Omega t)v(t)$. Under this transformation, 
  the equation (\ref{eq2}) becomes
    \begin{equation}\label{eq3}
  \frac{\partial V}{\partial t}+\nu AV =
       B_{\alpha}( \Omega\,t,V,V)+E_{\alpha}(\Omega\,t)f,
\end{equation}
where
\begin{eqnarray*}
  B_{\alpha}(\Omega\,t,V,V) &=&
    E_{\alpha}(\Omega\,t)P_L\{
    [ \mathcal{R}_{\alpha}E_{\alpha}(-\Omega\,t)V]\times
    [ \mbox{curl}(E_{\alpha}(-\Omega\,t)V)] \}\nonumber \\
   &=& -E_{\alpha}(\Omega\,t)B_{\alpha}( E_{\alpha}(-\Omega\,t)V,
        E_{\alpha}(-\Omega\,t)V).\label{trans16}
\end{eqnarray*}
For each wave number $n\in \mathbb{Z}^3$,  
\begin{equation*}\label{trans18}
  B_{\alpha}(\Omega\,t,V_k,V_m)_n =
    i\frac{1}{1+\alpha^2|\check{k}|^2}
    E_{\alpha}(\Omega\,t)_nP_n[E_{\alpha}(-\Omega\,t)_kV_k\times
     (\check{m}\times E_{\alpha}(-\Omega\,t)_mV_m)],
\end{equation*}
which is explicitly time-dependent with rapidly
varying
coefficients. This suggests that, for $\Omega >>1$, the dynamic mechanisms of (\ref{eq3})
evolve over 
two different time scales; the first one being
induced by the fast Poincar\'{e} waves and the second given by the evolution
of the Poincar\'{e} ``slow envelope'' $V(t)$.
$B_{\alpha}(\Omega\,t,V,V)$ contains resonant terms
($\Omega\,t$-independent terms) and nonresonant terms ($\Omega\,t$-dependent terms), and
we can decompose it as
\[
  B_{\alpha}(\Omega\,t,V,V)
   = \tilde{B}_{\alpha}(V,V) +
     B_{\alpha}^{osc}(\Omega\,t,V,V).
\]
Observe that  $B_{\alpha}^{osc}(\Omega\,t,V,V)$ contains all
nonresonant terms and $\tilde{B}_{\alpha}(V,V)$,
``the resonant bilinear operator'', contains all resonant
terms.
Averaging over fast time scale in the limit $\Omega\rightarrow\infty$ 
removes the nonresonant operator:
\[
  \lim_{\Omega\rightarrow\infty}\frac{1}{2\pi}
   \int_{0}^{2\pi}B_{\alpha}^{osc}(\Omega\,s,V,V)_n\, ds
    = 0,
\]
and we arrive at the {\em resonant
limit $\alpha$-equations};
\begin{eqnarray}
 \frac{\partial w}{\partial t}+\nu Aw =
    \tilde{B}_{\alpha}(w,w) +
    \tilde{f} && \label{eq4} \\
 w(0)=v(0)\hspace{.5in} && \label{decoms8}\nonumber
\end{eqnarray}
where
\begin{eqnarray*}
  \tilde{B}_{\alpha}(w,w)&=&\lim_{\Omega\rightarrow\infty}
    \frac{1}{2\pi}\int_{0}^{2\pi}B_{\alpha}(\Omega\,s,w,w)\, ds
   \label{decoms9}\\
  \tilde{f}&=&\lim_{\Omega\rightarrow\infty}\frac{1}{2\pi}
  \int_{0}^{2\pi}E_{\alpha}(\Omega\,s)f\, ds. \label{decoms10}
\end{eqnarray*}

The existence of regular solutions of the resonant limit $\alpha$-equations (\ref{eq4}) was 
established in \cite{BMN2, Kim2} based on rigorous a priori estimates of the (bilinear) resonant 
limit operator 
$\tilde{B}_{\alpha}(w,w)$. The estimates were uniformly in $\alpha$.
Bootstrapping from global regularity of the resonant limit  $\alpha$-equations,
the existence of a global regular solution of the full 3D RNS-$\alpha$ for large $\Omega$
(Theorem \ref{Th1} and \ref{Th2})
was proved. The convergence of the solutions to those of the exact 
RNS equations as $\alpha\rightarrow 0^+$ was also 
proved in the context of attractors \cite{Kim2}.

%%%%%%%%%%%%%%%%%%%%%%%%%%%%%%%%%%%%%%%%%%%%%%%%%%%%%%%%%%%%%%%%%%%

\section {Resonant Set and Operator Splitting} \label{operator-splitting}

%%%%%%%%%%%%%%%%%%%%%%%%%%%%%%%%%%%%%%%%%%%%%%%%%%%%%%%%%%%%%%%%%%%

Let $H$ be any Hilbert space and  $\overline{u}$  the $x_3$-averaging of $u\in H$:
\begin{equation}\label{baro1}
 \overline{u}(t,x_1,x_2) = \frac{1}{2\pi a_3}\int_0^{2\pi a_3}u(t,x_1,x_2,x_3)\,dx_3.\nonumber
\end{equation}
Denote $\overline{H}=\{ \overline{u}(t,x_1,x_2)| u\in H\}$. Then $\overline{H}$ is
a closed subspace of $H$, and any $u\in H$ has a unique
representation $u=\overline{u}+u^{\bot} \in \overline{H}\bigoplus H^{\bot}$. Note that
$\overline{u^{\bot}}=0$. 
This defines orthogonal projections
$P_b$ and $P_b^{\bot} = I-P_b$ on $H$ as
\begin{equation}\label{baro2}
  P_bu=\overline{u}\ \ \mbox{and}\ \  P_b^{\bot}u=u^{\bot} \nonumber
\end{equation}
We call $P_b$ a \emph{barotropic projection} and $P_b^{\bot}$ a \emph{baroclinic projection}, which 
make an orthogonal decomposition $H=\overline{H}\bigoplus H^{\bot}$ 
with $\overline{H}=P_bH$ and  $H^{\bot}=P_b^{\bot}H$ (For more details, see \S 3.4, \cite{Kim2}).

\begin{lemma}\label{lemma1}
 For any finite dimensional subspace $S\subset H$ with $dim(S) = d$, there exists an orthonormal 
 basis $\left\{ \phi_i\right\}_{1\leq i\leq d}$ in $H$, such that
 $\left\{ \phi_i\right\}_{1\leq i\leq d_1}$ and $\left\{ \phi_i\right\}_{d_1+1\leq i\leq d}$
 orthonormally span $P_bS$ and $P_b^{\bot}S$, respectively. (Lemma 4.2.2, \cite{Markus}).
\end{lemma}

Let's set $D_l(k,m,n)=\pm\omega_{\alpha k}\pm\omega_{\alpha m}\pm\omega_{\alpha n}$,
where $l=1,2,...,8$ is the combination of signs $\pm$.
The resonant nonlinear interactions of Poincar\'{e}
waves for $\tilde{B}_\alpha(w, w)$ in (\ref{eq4})
are present when the Poincar\'{e} frequencies satisfy the resonant
relation $D_l(k,m,n)=0$, and we define the corresponding \emph{resonant
set} $K$ by
\begin{equation*}
        K = \left\{ (k, m, n)\in \mathbb{Z}^3 :
   \pm\omega_{\alpha k}\pm\omega_{\alpha m}\pm\omega_{\alpha n} = 0,
    \ \ n=k+m\right\}.
\end{equation*}
We can decompose the resonant set $K$ into three groups for further analysis;
 pure 2D interactions ($K_{2D}$), two wave interactions ($\tilde{K}$),
 and three wave interactions ($K^*$). 
\begin{itemize}
  \item[(i)\ ] $K_{2D}=\{ (k,m,n)\in K|k_3=m_3=n_3=0\}$
    corresponds to pure two dimensional horizontal
    interactions
    (i.e., depends on $x_1, x_2$ and does not depend on $x_3$ in
    physical space.)
  \item[(ii)] $\tilde{K}=\{ (k,m,n)\in K |  k_3m_3n_3 = 0,
                  \ \ k_3^2+m_3^2+n_3^2 \neq 0\}$ is the set of two wave resonances.
         Here $ k_3m_3n_3 = 0$ represents that one or two of  $k_3, m_3$ and $ n_3$
         would be zero. But, if two of them are zero, we have 1-wave interaction
         which are excluded. It requires the second condition
         $k_3^2+m_3^2+n_3^2 \neq 0$. 
        This is the case when one of the three frequencies $\omega_\alpha$ equals zero
         and two remaining $\omega_\alpha$ are nonzero; for example,
         $\{ (k,m,n)\in K | \xi_{\alpha n}=0, \omega_{\alpha k}+\omega_{\alpha m}
         = 0, \omega_{\alpha k}\neq 0\neq \omega_{\alpha m}\} = K_{14}
         = (K_1\cap K_4)\setminus K_{2D}$.
         $\tilde{K}$ can be expressed in the way of
         \[ \tilde{K} = K_{14}\cup K_{24}\cup K_{34}, \]
         where $K_{j4}=(K_j\cap K_4)\setminus K_{2D}$ for $j=1,2,3$ and
        \begin{eqnarray*}
           K_{14} &=& \{ (k,m,n)\in K | n_3=0, \check{k}_3=-\check{m}_3\neq 0,
                 |\check{m}|=|\check{k}| \} \label{resoset8}\\
           K_{24} &=& \{ k_3=0, \check{m}_3=\check{n}_3\neq 0,
                     |\check{m}|=|\check{n}|\}\label{resoset9}\\
           K_{34} &=& \{ m_3=0, \check{k}_3=\check{n}_3\neq 0,
                     |\check{k}|=|\check{n}|\}\label{resoset10}
        \end{eqnarray*}
        Formally there exist three more 2-wave cones, but they are empty sets	\cite{BMN1}.
    \item[(iii)] $K^*=\{ (k,m,n)\in K | k_3m_3n_3 \neq 0 \}$ is
        the set of strict three wave resonances.
   \end{itemize}
 Then the resonant limit operator $\tilde{B}_{\alpha}(w,w)$ on $K$ has the following representation:
 \begin{equation*}
  \tilde{B}_{\alpha}(w,w) 
  = \tilde{B}_{I}^{\alpha}(\overline{w},\overline{w})
     + \tilde{B}_{II}^{\alpha}(\overline{w},w^{\bot})
     + \tilde{B}_{III}^{\alpha}(w^{\bot},w^{\bot})
\end{equation*}
where
\begin{eqnarray*}
  \tilde{B}_{I}^{\alpha}(\overline{w},\overline{w})
        &=& i\sum_{K_{2D}}\frac{1}{1+\alpha^2|\check{k}|^2}
       P_n[\overline{w}_k\times (\check{m}\times\overline{w}_m)]
        \label{baro12}\\
       &=&-i\sum_{K_{2D}}\frac{1}{1+\alpha^2|\check{k}|^2}
          P_n[(\check{m}\cdot\overline{w}_k)\overline{w}_m]\\
  \tilde{B}_{II}^{\alpha}(\overline{w},w^{\bot})&=&
         -i\sum_{K_{24}}\frac{1}{1+\alpha^2|\check{k}|^2}
         P_n[(\check{m}\cdot\overline{w}_k)w_m^{\bot}] \\
          & & +\frac{i}{2}\sum_{K_{24}}\frac{1}{1+\alpha^2|\check{k}|^2}
           P_n[(\overline{w}_k\cdot w_m^{\bot})\check{m}+(\check{k}\cdot w_m^{\bot})\overline{w}_k] \\
         & & +\frac{i}{2}\sum_{K_{34}}\frac{1}{1+\alpha^2|\check{n}|^2}P_n
            [(w_k^{\bot}\cdot\overline{w}_m)\check{m}-(\check{m}\cdot w_{k}^{\bot})\overline{w}_m ] \\
         & & -\frac{i}{2}\sum_{K_{34}}\frac{1}{1+\alpha^2|\check{n}|^2}
          \frac{1}{|\check{n}|^2}P_n[\check{n}\times
           \left\{ (\check{k}\times w_k^{\bot})\times
             (\check{m}\times \overline{w}_m)\right\} ] \label{baro13} \\
  \tilde{B}_{III}^{\alpha}(w^{\bot},w^{\bot})&=&
    \sum_{(k,m,n)\in K^{*}}Q_{kmn}(w_k^{\bot},w_m^{\bot}),
\end{eqnarray*}
where $Q_{kmn}(v_k,v_m)$
is a bilinear form in $v_k,v_m\in\mathbb{C}^3$
with the estimate
\[
   |Q_{kmn}(v_k,v_m)|
    \leq |\check{m}||v_k||v_m|.
\]

In this paper, we consider the catalytic $\alpha$-limit system, that is, the resonant $\alpha$-limit equations 
 not including the strict three-wave resonant operator $B_{III}^\alpha$. As pointed out in \cite{BMN1}, there exist a
 generic set of domain parameters $(a_1, a_2, a_3)$, dense in $\mathbf{R}_3^+$ and of full 
 Lebesque measure, for which $B_{III}^\alpha$ is identically zero. 
 Defining 
 \[  \tilde{B}_c^{\alpha}(w,w) = B_I^{\alpha}(\bar{w},\bar{w})+B_{II}^{\alpha}(\bar{w},w^\bot ) \]
 the resonant limit $\alpha$-equations (\ref{eq4}) take the form 
\begin{eqnarray}
 \frac{dw}{dt}+\nu Aw+\tilde{B}_c^{\alpha}(w,w) = \tilde{f} & &  \label{eq5}\\ 
  w(0) =w_0,\ \ \ \nabla\cdot w=0, \nonumber & &
\end{eqnarray}
where $\tilde{B}_c^{\alpha}$ is a bilinear operator of the catalytic system;
\begin{eqnarray*}
   B_I^{\alpha}(\bar{w},\bar{w})&=&P_L(\mathcal{R}_{\alpha}\bar{w}\cdot\nabla_h)\bar{w}
     = P_L[\mathcal{R}_{\alpha}\overline{w}\times\mathrm{curl}\,\overline{w}]\\
   B_{II}^{\alpha}(\bar{w},w^\bot ) &=&  P_L(\mathcal{R}_{\alpha}\bar{w}\cdot\nabla) w^\bot
      = P_L[\mathcal{R}_{\alpha}\overline{w}\times\mathrm{curl}\,w^\bot].
\end{eqnarray*}

Eqs.(\ref{eq5}) has a unique solution 
$w\in L^\infty(0, T;H^1)\cap L^2(0, T; H^2)$ for $\alpha\geq 0$ with 
$T<\infty$ \cite{Kim2, Markus}. Hence the semigroup $S_\alpha (t): H\rightarrow H$ is defined:
$S_\alpha (t)w_0 = w(t)$, where $w(t)$ is the solution of (\ref{eq5}). 
The semigroup $S_\alpha (t)$ has an absorbing ball in $H$ and a global attractor 
$\mathcal{A}_\alpha \subset H$ \cite{Kim2}.

%%%%%%%%%%%%%%%%%%%%%%%%%%%%%%%%%%%%%%%%%%%%%%%%%%%%%%%%%%%%%%%%%%%

\section {Attractor Dimensions}

%%%%%%%%%%%%%%%%%%%%%%%%%%%%%%%%%%%%%%%%%%%%%%%%%%%%%%%%%%%%%%%%%%%

We now estimate the dimension of the global attractor $\mathcal{A}_\alpha$.
Attractor dimension is associated with the number of degrees of freedom for turbulent rotating fluids
for the Coriolis parameter $\Omega$ large.
The focus is on the estimation of the dimensions uniform in $\alpha$, no blow up as $\alpha\rightarrow 0^+$.

We consider the variational equation corresponding to Eq.\,(\ref{eq5}):
\begin{eqnarray}
 \frac{d\Phi}{dt} &=& L_c(w)\Phi \label{var_eq1} \\
     \Phi(0) &=& \xi \nonumber
\end{eqnarray}
 with
\begin{eqnarray*}
L_c(w)\Phi &=& -\nu A\Phi-\tilde{B}_c^\alpha (w,\Phi)-
     \tilde{B}_c^\alpha (\Phi,w) \\
   &=& -\nu A\Phi -B_I^\alpha (\bar{w},\bar{\Phi})-B_I^\alpha (\bar{\Phi},\bar{w})
   -B_{II}^\alpha (\bar{w},\Phi^\bot )-B_{II}^\alpha (\bar{\Phi},w^\bot )
\end{eqnarray*}
Following the standard procedure as in section 13.4, \cite{Robinson1}, we can show that the equations
(\ref{var_eq1}) have a unique solution:
\begin{lemma}\label{lemma1}
 If $w$ is a solution of Eq.\,(\ref{eq5}), then Eq.\,(\ref{var_eq1})
 possesses a unique solution
 \[ \Lambda (t,w_0)\xi = \Phi (t) \in L^2(0, T; H^1)\cap C([0, T]; H),\ \forall T>0. 
  \]
 Furthermore, for every $t>0$, the flow $w_0\rightarrow S_{\alpha t}(w_0)$ 
 generated by the limit $\alpha$-equations (\ref{eq5})
   is uniformly differentiable on $\mathcal{A}_\alpha$ with the differential 
   \[ DS_{\alpha t}(w_0) = \Lambda (t, w_0): \xi\in \mathcal{A}_\alpha 
      \rightarrow \Phi(t)\in H, \]
   where $\Phi$ is the solution of (\ref{var_eq1}). In particular, the linear operator
   $\Lambda(t, \bar{w}_0)$ is compact for all $t>0$.
\end{lemma}
%{\bf Proof.} See Appendix \ref{app2}.

Now we are ready to  apply the following theorem to estimate the fractal dimension of the global attractor.
Set 
\[ \mathcal{TR}_n(\mathcal{A}_\alpha) = \sup_{w_0\in\mathcal{A}_\alpha}\sup_{P^{(n)}(0)}
    \limsup_{t\rightarrow\infty} \frac{1}{t}\int_0^t\mathrm{Tr}\left(L(s;w_0)P^{(n)}(s)\right)\,ds, \]
where $\mathrm{Tr}(M)$ is the trace of $M$.
\begin{theorem}\label{theorem1}(\cite{Robinson1}, Theorem 13.16)
  Suppose that $S_{\alpha t}(t)$ is uniformly differentiable on $\mathcal{A}_\alpha$ and that there 
  exists a $t_0$ such that $\Lambda (t, w_0)$ is compact for all $t\geq t_0$.
  If $\mathcal{TR}_n(\mathcal{A}_\alpha ) < 0$ then $\mathrm{dim_f}(\mathcal{A}_\alpha )\leq n$. 
\end{theorem}

Let $\Phi_1, ..., \Phi_N$ be solutions of the linearized system (\ref{var_eq1})
with corresponding initial conditions $\xi_1, ..., \xi_N$.
Let $\phi_1, ..., \phi_N$ be the orthonormal system spanning
$\{ \Phi_1, ..., \Phi_N\}$. Then the trace at time $t  \geq 0$ can be sought like
\begin{equation}\label{trace}
Tr\left( L_c(t)P_N(\Phi_1(t), ..., \Phi_N(t))\right)
  = \sum_{i=1}^{N}\left< L_c(t)\phi_i(t), \phi_i(t)\right>
\end{equation}

\noindent Using $\phi = \bar{\phi}+\phi^\bot$ and bilinearity of the operators $B_I^\alpha$ and
$B_{II}^\alpha$, the inner product on the right-handed side of (\ref{trace}) can be expressed as
\begin{eqnarray*}
-\left< L_c(t)\phi_i, \phi_i\right> 
       &=& \left< \nu A\bar{\phi}_i+
         B_I^\alpha (\bar{w},\bar{\phi}_i)+B_I^\alpha (\bar{\phi}_i,\bar{w}),
         \bar{\phi}_i\right> \\
      & & \ \   +\left< \nu A\phi_i^\bot +
         B_{II}^\alpha (\bar{w},\phi_i^\bot )+B_{II}^\alpha (\bar{\phi}_i,w^\bot ),
         \bar{\phi}_i\right> \\
       & & \ \ + \left< \nu A\bar{\phi}_i+
         B_I^\alpha (\bar{w},\bar{\phi}_i)+B_I^\alpha (\bar{\phi}_i,\bar{w}),
         \phi_i^\bot \right> \\
      & & \ \  +\left< \nu A\phi_i^\bot +
         B_{II}^\alpha (\bar{w},\phi_i^\bot )+B_{II}^\alpha (\bar{\phi}_i,w^\bot ),
         \phi_i^\bot \right> 
\end{eqnarray*}

$\tilde{B}_\alpha$ is the skew-symmetric limit bilinear form. We collect the following properties
from \cite{BMN1, BMN2, Kim2, Markus}. 

\begin{lemma}\label{bilinear-property} For any $w, \phi \in H^1$
 \begin{enumerate}
   \item $\left< \nu A\bar{\phi}_i, \phi_i^\bot \right> =
          \nu\left< A\bar{\phi}_i, \phi_i^\bot \right> =
          \nu\left< P_b(A\bar{\phi}_i), P_b^\bot (A\phi_i^\bot )\right> = 0$.
   \item $\left< \nu A\phi_i^\bot , \bar{\phi_i} \right> = 0$. Similar to 1.
   \item $B_I^\alpha (\bar{w},\bar{\phi}_i) = P_b\tilde{B}_\alpha (w,\phi_i)$.
   \item $\left< B_I^\alpha (\bar{w},\bar{\phi}_i), \phi_i^\bot\right> =
        \left< B_I^\alpha (\bar{w},\bar{\phi}_i), P_b^\bot\phi_i\right> =
        \left< P_b^\bot B_I^\alpha (\bar{w},\bar{\phi}_i), \phi_i\right> = \left<0,
        \phi_i\right>= 0$, since
        $P_b$ and $P_b^\bot$ are self-adjoint.
   \item Similarly, $\left< B_I^\alpha (\bar{\phi}_i, \bar{w}), \phi_i^\bot\right> = 0$.
 \end{enumerate}
\end{lemma}
\noindent With Lemma \ref{bilinear-property} we
can reduce the above inner product.
\begin{eqnarray*}
 -\left< L_c(t)\vec{\phi}_i , \vec{\phi}_i\right>
  & = &\left< \nu A\bar{\phi}_i+
         B_I^\alpha (\bar{w},\bar{\phi}_i)+B_I^\alpha (\bar{\phi}_i,\bar{w}),
         \bar{\phi}_i\right> \\
      & & \ \ \    + \left< B_{II}^\alpha (\bar{w},\phi_i^\bot )
            +B_{II}^\alpha (\bar{\phi}_i,w^\bot ),
         \bar{\phi}_i\right> \\
      & & \ \ \ 
          +\left< \nu A\phi_i^\bot +
          B_{II}^\alpha (\bar{w},\phi_i^\bot )+B_{II}^\alpha (\bar{\phi}_i,w^\bot ),
         \phi_i^\bot \right>
\end{eqnarray*}

\begin{corollary} (\cite{Markus}, Cor 3.2.13, p45)
If $f\in P_bH, h\in \bar{H}$ and $g\in P_b^\bot H\cap V$, then
  \[ \left< B_{II}^\alpha (f,g), h\right> = 0. \]
\end{corollary}

\textbf{Proof.} From Lemma 3.4 in \cite{Kim2}
   with the self-adjoint property of $P_b$.
   ($h\in \bar{H}$ so that $h=P_bh$). \ \ $\blacksquare$
\medskip

From this corollary we obtain
\[ \left< B_{II}^\alpha (\bar{w}, \phi_i^\bot )+
  B_{II}^\alpha (\bar{\phi}_i, w^\bot ), \bar{\phi}_i\right> = 0.
\]
From Theorem 4.2 of \cite{Kim2}
\[ \left< B_{II}^\alpha (\bar{w}_i, \phi_i^\bot ), \phi_i^\bot\right> =0.
\]
Thus, by the construction of $(\phi_i)_{1\leq i\leq N}$, for each $i$ 
either $\phi_i = 0$ or $\phi_i^\bot = 0$ 
\[ \left< B_{II}^\alpha (\bar{\phi}_i, w^\bot ), \phi_i^\bot\right> =0.
\]
(\cite{Markus}, p58 with Lemma 4.2.2),
and
\begin{equation*}
 -\left< L_c(t)\vec{\phi}_i, \vec{\phi}_i\right> =
   \left< \nu A\bar{\phi}_i+B_I^\alpha (\bar{w},\bar{\phi})+
         B_I^\alpha (\bar{\phi}_i,\bar{w}), \bar{\phi}_i\right>
         +\left< \nu A\phi_i^\bot, \phi_i^\bot\right> .
\end{equation*}
 In summary,
\begin{lemma}\label{trace-lemma}
  Let $\Phi_1, ..., \Phi_N$ be solutions of the linearized system (\ref{var_eq1})
with corresponding initial conditions $\xi_1, ..., \xi_N$.
Let $\phi_1, ..., \phi_N$ be the orthonormal system spanning
$\{\Phi_1, ..., \Phi_N\}$ in $H$. Then the trace at time $t  \geq 0$ is given by
\begin{eqnarray*}
Tr\left( L_c(t)P_N(\Phi_1(t), ..., \Phi_N(t))\right)
  = \sum_{i=1}^N \left< L_c(t)\phi_i, \phi_i\right>
    \hspace{1in}  & & \\
  =  - \sum_{i=1}^N \left[ \left< \nu A\bar{\phi}_i+B_I^\alpha (\bar{w},\bar{\phi})+
         B_I^\alpha (\bar{\phi}_i,\bar{w}), \bar{\phi}_i\right>
         +\left< \nu A\phi_i^\bot, \phi_i^\bot\right> \right] & &\\
\end{eqnarray*}
\end{lemma}

Now let's estimate for the Time-Averaged Trace in Lemma \ref{trace-lemma}.

%---------------------------
 \subsection{A Priori Estimates}
%---------------------------

%---------------------------
  \subsubsection{Estimate 1}
%---------------------------

Refer to Theorem 4.2.4 (p58) in \cite{Markus}. We may use
   \[ \sum_{i=1}^{N_1} \left< B_I^\alpha (\bar{\phi}_i,\bar{w}),\bar{\phi}_i\right>
      = \int_{T^2}\sum_{i=1}^{N_1}\left( \mathcal{R}_\alpha\bar{\phi}_i\times
         \mathrm{curl} \bar{w}\right)\cdot\bar{\phi}_i\,dx.
   \]
Instead, we will follow the line below. 
Denote $\bar{\phi}_i = <\bar{\phi}_{i1}, \bar{\phi}_{i2}, \bar{\phi}_{i3}>$.
By Lemma \ref{bilinear-relation},
  \begin{eqnarray*}
     \sum_{i=1}^{N_1} \left< B_I^\alpha (\bar{\phi}_i,\bar{w}),\bar{\phi}_i\right>
       = \sum_{i=1}^{N_1} \left[ \left< B(\mathcal{R}_\alpha\bar{\phi}_i, \bar{w}), \bar{\phi}_i\right>
       - \left< B(\bar{\phi}_i,\bar{w}), \mathcal{R}_\alpha\bar{\phi}_i\right>\right] 
       \hspace{1in} &&\\
      = \sum_{i=1}^{N_1} \left[\left< P_L\left[\left(\mathcal{R}_\alpha\bar{\phi}_i
       \cdot\nabla\right)\bar{w}\right], \bar{\phi}_i\right> -
       \left< P_L\left[\left(\bar{\phi}_i
       \cdot\nabla\right)\bar{w}\right], \mathcal{R}_\alpha \bar{\phi}_i\right>\right]
       \hspace{1.2in} &&\\
     = \sum_{i=1}^{N_1} \left[\left< \left(\mathcal{R}_\alpha\bar{\phi}_i
       \cdot\nabla\right)\bar{w}, \bar{\phi}_i\right> -
       \left< \left(\bar{\phi}_i
       \cdot\nabla\right)\bar{w}, \mathcal{R}_\alpha \bar{\phi}_i\right>\right]
       \hspace{1.9in}  & &\\
         \ \ \text{(since $P_L$ is self-adjoint and $P_L\bar{\phi}_i = \bar{\phi}_i$.)} 
         \hspace{.7in} &&\\
     = \sum_{i=1}^{N_1} \left[ \int_{T^2}\sum_{j=1}^2\sum_{k=1}^3\left( \mathcal{R}_\alpha\bar{\phi}_{i}\right)_j
       D_j\bar{w}_k \bar{\phi}_{ik}\,dx -
       \int_{T^2}\sum_{j=1}^2\sum_{k=1}^3\bar{\phi}_{ij}
       D_j\bar{w}_k \left( \mathcal{R}_\alpha \bar{\phi}_{i}\right)_k\,dx\right] &&\\
     = \int_{T^2}\sum_{i=1}^{N_1}\sum_{j=1}^2\sum_{k=1}^3
       \left[ \left( \mathcal{R}_\alpha\bar{\phi}_{i}\right)_j
        D_j\bar{w}_k \bar{\phi}_{ik} - \bar{\phi}_{ij}
       D_j\bar{w}_k \left( \mathcal{R}_\alpha \bar{\phi}_{i}\right)_k\right]\,dx
       \hspace{1in} &&
  \end{eqnarray*}
Then,
{\small
\begin{eqnarray*}
 \left|\int_{T^2}\sum_{i=1}^{N_1}\sum_{j=1}^2\sum_{k=1}^3
       \left[ \left( \mathcal{R}_\alpha\bar{\phi}_{i}\right)_j
        D_j\bar{w}_k \bar{\phi}_{ik} - \bar{\phi}_{ij}
       D_j\bar{w}_k \left( \mathcal{R}_\alpha \bar{\phi}_{i}\right)_k\right]\,dx\right|
       \leq \hspace{3in} & &  \nonumber \\
     \hspace{0.5in}   \left| \int_{T^2}\sum_{i=1}^{N_1}\sum_{j=1}^2\sum_{k=1}^3
       \left[ \left( \mathcal{R}_\alpha\bar{\phi}_{i}\right)_j
        D_j\bar{w}_k \bar{\phi}_{ik} \right]\,dx\right|  %\nonumber \\
          + \left| \int_{T^2}\sum_{i=1}^{N_1}\sum_{j=1}^2\sum_{k=1}^3
       \left[  \bar{\phi}_{ij}
       D_j\bar{w}_k \left( \mathcal{R}_\alpha \bar{\phi}_{i}\right)_k\right]\,dx\right|  \hspace{1.7in} \label{eq4.1}
\end{eqnarray*}
}

We only need to estimate the first term on the right-handed side because the second term
can be estimates similarly.
\begin{equation*}
  \left| \int_{T^2}\sum_{i=1}^{N_1}\sum_{j=1}^2\sum_{k=1}^3
       \left[ \left( \mathcal{R}_\alpha\bar{\phi}_{i}\right)_j
        D_j\bar{w}_k \bar{\phi}_{ik} \right]\,dx\right|
        \leq
         \int_{T^2}\left| \sum_{i=1}^{N_1}\sum_{j=1}^2\sum_{k=1}^3
       \left[ \left( \mathcal{R}_\alpha\bar{\phi}_{i}\right)_j
        D_j\bar{w}_k \bar{\phi}_{ik} \right]\right|\,dx
\end{equation*}
Here,
\begin{eqnarray*}
  \left| \sum_{i=1}^{N_1}\sum_{j=1}^2\sum_{k=1}^3
       \left[ \left( \mathcal{R}_\alpha\bar{\phi}_{i}\right)_j
        D_j\bar{w}_k \bar{\phi}_{ik} \right]\right|
  \leq
    \sum_{i=1}^{N_1}\left| \sum_{j=1}^2\sum_{k=1}^3
       \left[ \left( \mathcal{R}_\alpha\bar{\phi}_{i}\right)_j
        D_j\bar{w}_k \bar{\phi}_{ik} \right]\right| \hspace{1in} & &\\
   \leq
    \sum_{i=1}^{N_1}\left\{\left[\sum_{j=1}^2\sum_{k=1}^3 \left( D_j\bar{w}_k(x)\right)^2\right]^{1/2}
      \left[ \sum_{j=1}^2\sum_{k=1}^3 \left(\mathcal{R}_\alpha\bar{\phi}_{i}(x)\right)_j^2
      \bar{\phi}_{ik}^2(x) \right]^{1/2}\right\} \hspace{.6in} &&  \\
  = \left|\nabla\bar{w}(x)\right| \sum_{i=1}^{N_1} \left[
     \sum_{j=1}^2\left(\mathcal{R}_\alpha\bar{\phi}_{i}(x)\right)_j^2
     \sum_{k=1}^3 \bar{\phi}_{ik}^2(x) \right]^{1/2} \hspace{1.9in} && \\
  \leq \left|\nabla \bar{w}(x)\right| \sum_{i=1}^{N_1}
    \left[ \sum_{j=1}^3\left(\mathcal{R}_\alpha\bar{\phi}_{i}(x)\right)_j^2
        \sum_{k=1}^3\bar{\phi}_{ik}^2(x)\right]^{1/2} \hspace{1.9in} && \\
  = \left|\nabla\bar{w}(x)\right| \sum_{i=1}^{N_1} \left\{ \left[
     \sum_{j=1}^3\left(\mathcal{R}_\alpha\bar{\phi}_{i}(x)\right)_j^2\right]^{1/2}
     \left[\sum_{k=1}^3 \bar{\phi}_{ik}^2(x) \right]^{1/2} \right\} \hspace{1.3in} &&\\
 \leq \left|\nabla\bar{w}(x)\right| \left[ \sum_{i=1}^{N_1}
     \sum_{j=1}^3\left(\mathcal{R}_\alpha\bar{\phi}_{i}(x)\right)_j^2\right]^{1/2}
     \left[\sum_{i=1}^{N_1}\sum_{k=1}^3 \bar{\phi}_{ik}^2(x) \right]^{1/2} \hspace{1.3in} && \\
  = \left|\nabla\bar{w}(x)\right| \left[ \sum_{i=1}^{N_1}
     \left|\mathcal{R}_\alpha\bar{\phi}_{i}(x)\right|^2\right]^{1/2}
     \left[\sum_{i=1}^{N_1}\left|\bar{\phi}_{i}(x)\right|^2 \right]^{1/2} \hspace{1.7in} &&  \\
   \leq \left|\nabla\bar{w}(x)\right| \left[ \sum_{i=1}^{N_1}
     \left|\mathcal{R}_\alpha\bar{\phi}_{i}(x)\right|^2\right]^{1/2}
     \left(\sum_{i=1}^{N_1}\left|\bar{\phi}_{i}(x)\right|^2 \right) \hspace{1.7in} &&
\end{eqnarray*}
(Why do we need the last inequality? It is to get the estimate $||\rho||_0^2\leq c_l\sum||\bar{\phi}_i||^2$
 so that we can combine two nonlinear-term estimates to have $||\bar{w}||^2$ and be able to estimate
 $\int_t^{t+\tau}||\bar{w}(s)||^2\,ds$. Otherwise, we will have $||\rho^{1/2}||_0 =  N_1$ and get additional
 $||\bar{w}||$ term and have to estimate  $\int_t^{t+\tau}||\bar{w}(s)||\,ds$:
 $ ||\rho^{1/2}||_0^2 = \int_\Omega|\rho^{1/2}(x)|^2\,dx = 
     \sum_{i=1}^{N_1}\int_\Omega |\bar{\phi}_i(x)|^2dx = \sum_{i=1}^{N_1}||\bar{\phi}_i||_0^2 = N_1
 $).
 \smallskip
 
\noindent Integrating both sides and using H\"{o}lder's inequality, we obtain
{\small 
\begin{eqnarray*}
 \int_{T^2} |\nabla\bar{w}(x)| \left[ \sum_{i=1}^{N_1}
     |\mathcal{R}_\alpha\bar{\phi}_{i}(x)|^2\right]^{1/2}
     \left[\sum_{i=1}^{N_1}|\bar{\phi}_{i}(x)|^2 \right]\,dx
     &\leq&
     ||\rho_\alpha||_{\infty} \int_{T^2} |\nabla\bar{w}(x)|\rho(x)\,dx\\
     &\leq& ||\rho_\alpha ||_{\infty} ||\bar{w}|| |\rho|_{L^2}, 
\end{eqnarray*}
}
 where $\rho_\alpha (x)= \left[ \sum_{i=1}^{N_1}|\mathcal{R}_\alpha\bar{\phi}_i(x)|^2\right]^{1/2}$
  and $\rho (x)= \sum_{i=1}^{N_1}|\bar{\phi}_i(x)|^2$ with
   $||\rho_\alpha||_\infty = \sup_{x\in \mathbb{T}^3}|\rho_\alpha(x)|$.
Therefore, we get the estimate
\begin{equation}\label{eq4.4}
  \left|\sum_{i=1}^{N_1}\left< B_I^\alpha (\bar{\phi}_i, \bar{w}), \bar{\phi}_i\right>\right|
  \leq
  2||\rho_\alpha||_{\infty} |\rho|_{L^2}||\bar{w}||  
\end{equation}

%--------------------------
  \subsubsection{Estimate 2 ($||\rho_\alpha||_\infty$ Estimate)} 
%---------------------------

Let $\theta_i = \overline{\phi}_i$ and $v_i = (1-\alpha^2\Delta )^{-1}\theta_i$. Then,
\[ \rho_\alpha (x)  = \left[ \sum_{i=1}^{N_1}\left| (1-\alpha^2\Delta )^{-1}\theta_i\right|^2\right]^{1/2} 
     = \left[ \sum_{i=1}^{N_1}\left| v_i\right|^2\right]^{1/2} 
\]

For $\theta = \theta_i\in H$, we have $v = v_i\in H^3$. With the Sobolev embedding theorem, we can infer that
the existence of a dimensionless constant $c(\alpha)$ depending on $\alpha$ such that  
\[ 
  ||v||_\infty \leq C||v||_2 = C||(1-\alpha^2\Delta )^{-1}\theta||_2 \leq c(\alpha)||\theta||_0.
\]
Now we will compute the conatnat $c(\alpha)$, following the same process as in \cite{Ilyin1}.
Suppose that $\xi_1, \xi_2, ..., \xi_m\in \mathbb{R}$ and $\sum_{j=1}^m\xi_j^2 = 1$. 
Then, using the orthonormality of the $\theta_j$ $(j=1, 2, ..., m=N_1)$, and the above inequality we obtain
\begin{eqnarray*}
 \left|\sum_{j=1}^m\xi_jv_j(x)\right| &\leq & c(\alpha) ||\sum_{j=1}^m\xi_j\theta_j|| \\
     &=& c(\alpha) \left[ \int \left( \sum_{j=1}^m\xi_j\theta_j(x)\right)^2\,dx\right]^{1/2} \\
     &=& c(\alpha) \left(\sum_{i,j=1}^m\xi_i\xi_j\delta_{ij}\right)^{1/2} 
       =  c(\alpha) \left(\sum_{j=1}^m \xi_j^2\right)^{1/2} = c(\alpha)
\end{eqnarray*}

Using the representation $v_j = v_{j1}\cdot e_1 + v_{j2}\cdot e_2$ we find that 
\[ 
 \left[ \sum_{j=1}^m\xi_jv_j(x)\right]^2 \leq \left(\sum_{j=1}^m\xi_jv_{j1}(x)\right)^2 + 
     \left(\sum_{j=1}^m\xi_jv_{j2}(x)\right)^2 \leq c(\alpha)^2
\]
First, we set $\xi_j = v_{j1}(x)/(\sum_{j=1}^m(v_{j1}(x))^2)^{1/2}$ and later set  
$\xi_j = v_{j2}(x)/(\sum_{j=1}^m(v_{j2}(x))^2)^{1/2}$. 
Substituting these into the above inequality one after another,
we obtain
\[
  \rho_\alpha^2(x) = \sum_{j=1}^m |v_j(x)|^2 = \sum_{j=1}^m (v_{j1}(x))^2 + \sum_{j=1}^m (v_{j2}(x))^2
  \leq 2c(\alpha )^2.
\]

To compute $c(\alpha )$, we use the Fourier Series
\[ \theta_j (x) = \sum_{k\in Z_0^2} a_{jk}e^{ik\cdot x},\ \ \  Z_0^2 = \mathbb{Z}^2\{ 0\},\ 
  \ x = (x_1, x_2), \ \ k = (k_1, k_2) \]
so that 
\[ (1-\alpha^2\Delta )^{-1}\theta_j (x) = \sum_{k\in Z_0^2}\frac{a_{jk}}{{1+\alpha^2|k|^2}} e^{ik\cdot x}. \]
Thus,
\begin{eqnarray*}
  |(1-\alpha^2\Delta )^{-1}\theta_j (x)| &\leq & \sum_{k\in Z_0^2}\frac{|a_{jk}|}{1+\alpha^2|k|^2} \\
    &\leq &  \left( \sum_{k\in Z_0^2}\frac{1}{(1+\alpha^2|k|^2)^2}\right)^{1/2} 
            \left( \sum_{k\in Z_0^2} |a_{jk}|^2\right)^{1/2} \\
          &=& c_j(\alpha ) ||\theta_j||_0 \\
          &=& c_j(\alpha ), \ \ \ \ (\mbox{since}\ ||\theta_j||_0=1),
\end{eqnarray*}
where 
\begin{eqnarray*}
 c_j^2(\alpha ) &=&  \sum_{k\in Z_0^2}\frac{1}{(1+\alpha^2|k|^2)^2}\\
       &=&  \sum_{p=1}^\infty \frac{1}{(1+\alpha^2\lambda_p)^2},\\
       & & \ \ \{ \lambda_p, p=1, 2, ...\} = 
           \{ \hat{k}^2 = \hat{k}_1^2 + \hat{k}_2^2, \hat{k}_i = k_i/a_i \ \mbox{for}\ i=1,2, (k_1, k_2)\in Z_0^2\}\\
       &\leq & \sum_{p=1}^\infty \frac{1}{(1+\alpha^2c_1p)^2} \ \ \mbox{for an absolute constant}\ c_1\\
       &=& \frac{1}{(1+\alpha^2c_1)^2} + \sum_{p=2}^\infty \frac{1}{(1+\alpha^2c_1p)^2} \\
       &\leq & \frac{1}{(1+\alpha^2c_1)^2} + \int_1^\infty \frac{dx}{(1+\alpha^2c_1x)^2} \\
       &=& \frac{1}{1+\alpha^2c_1}\left[ \frac{1}{1+\alpha^2c_1} + \frac{1}{\alpha^2c_1}\right]  
\end{eqnarray*}

Then we can set
\[  \sum_{j=1}^{N_1}c_j^2(\alpha) \leq \frac{N_1}{1+\alpha^2c_1}\left[
   \frac{1}{1+\alpha^2c_1}+\frac{1}{\alpha^2c_1}\right] .
\]
Therefore, 
\[ ||\rho_\alpha||_{\infty}^2 \leq 2N_1 c^2(\alpha), \]
where 
\[
   c^2(\alpha) = \frac{1}{1+\alpha^2c_1}\left[\frac{1}{1+\alpha^2c_1}+\frac{1}{\alpha^2c_1}\right] .
\]
Observe that $\lim_{\alpha\rightarrow 0^+} c^2(\alpha) = 2 <\infty$. 

%----------------------
  \subsubsection{Estimate 3}
%----------------------

This is a new term that the exact Navier-Stokes equations don't have.
\begin{eqnarray*}
  \sum_{i-1}^{N_1}\left< B_I^\alpha (\bar{w},\bar{\phi}_i, \bar{\phi}_i\right>
     &=& \sum_{i=1}^{N_1}\left[ \left<B(\mathcal{R}_\alpha\bar{w},\bar{\phi}_i),
       \bar{\phi}_i\right>
      - \left< B(\bar{\phi}_i, \bar{\phi}_i),
  \mathcal{R}_\alpha\bar{w}\right>\right]\\
    & & \left( \text{Note that}\
      \left<B(\mathcal{R}_\alpha\bar{w},\bar{\phi}_i),
       \bar{\phi}_i\right> = 0 \ \text{by skew-symmetry}\right)\\
    &=& \sum_{i=1}^{N_1}\left<B(\bar{\phi}_i, \mathcal{R}_\alpha\bar{w}),
       \bar{\phi}_i\right> \ \ \text{by skew-symmetry}\\
    &=& \sum_{i=1}^{N_1}\left[ \int_{T^2} \sum_{j=1}^2\sum_{k=1}^3
      \bar{\phi}_{ij}D_j (\mathcal{R}_\alpha\bar{w}_k)\bar{\phi}_{ik}\,dx\right]
\end{eqnarray*}
Then,
\begin{eqnarray*}
 \left|\sum_{i=1}^{N_1}\left< B_I^\alpha (\bar{w},\bar{\phi}_i, \bar{\phi}_i\right>\right|
   &=& \left|
     \int_{T^2} \sum_{i=1}^{N_1}\sum_{j=1}^2\sum_{k=1}^3
      \bar{\phi}_{ij}D_j (\mathcal{R}_\alpha\bar{w}_k)\bar{\phi}_{ik}\,dx\right|\\
   &\leq& \int_{T^2}\left|\nabla (\mathcal{R}_\alpha\bar{w}_k)\right|
     \left[ \sum_{j=1}^2  \sum_{k=1}^3
       \left( \sum_{i=1}^{N_1}\bar{\phi}_{ij}(x)
       \bar{\phi}_{ik}^2(x)\right)^2\right]^{1/2}\,dx  \\
   &\leq& \int_{T^2}\left|\nabla (\mathcal{R}_\alpha\bar{w}_k)\right|
     \left[ \sum_{j=1}^3  \sum_{k=1}^3
       \left( \sum_{i=1}^{N_1}\bar{\phi}_{ij}(x)
       \bar{\phi}_{ik}^2(x)\right)^2\right]^{1/2}\,dx  \\
   &\leq& \int_{T^2} \left|\nabla (\mathcal{R}_\alpha\bar{w}_k)\right|
     \rho (x)\,dx, \ \ \ \ \ \text{where}\ \
     \rho (x)= \sum_{i=1}^{N_1} |\bar{\phi}_{i}(x)|^2\\
   &\leq& ||\mathcal{R}_\alpha\bar{w}||\,|\rho|_{L^2} \\
   &\leq& ||\bar{w}||\,|\rho|_{L^2}\\ 
   &\leq& ||\bar{w}||\left( c_l\sum_{i=1}^{N_1} ||\bar{\phi}_i||^2\right)^{1/2}\\
  & & \ \ \  \text{(by the Lieb-Thirring inequality; see p59 \cite{Markus} for detail),}\\
\end{eqnarray*}
where $c_l$ is an absolute constant and
\[ \left|\nabla (\mathcal{R}_\alpha\bar{w}_k)\right| =
    \left( \sum_{j=1}^2\sum_{k=1}^3\left|D_j(\mathcal{R}_\alpha\bar{w}_k)\right|\right)^{1/2}.
\]

%------------------------------
  \subsubsection{Collection of Estimate 1, 2, and 3}
%-------------------------------

\begin{eqnarray*}
  Tr (L_1(t)P_N(\bar{\Phi}_1(t),...,\bar{\Phi}_N(t)) \hspace{5in} &&\\ 
   \leq
    -\nu \sum_{i=1}^{N_1}||\bar{\phi}_i(t)||^2+ 
     2||\rho_\alpha||_{\infty}\left( c_l\sum_{i=1}^{N_1}||\bar{\phi}_i||^2\right)^{\frac{1}{2}}||\bar{w}|| 
     +||\bar{w}||\left( c_l\sum_{i=1}^{N_1} ||\bar{\phi}_i||^2\right)^{\frac{1}{2}} \hspace{1.5in}  && \\
  \leq -\nu \sum_{i=1}^{N_1}||\bar{\phi}_i(t)||^2+2\sqrt{2N_1}c(\alpha)||\bar{w}||
    \left( c_l\sum_{i=1}^{N_1} ||\bar{\phi}_i||^2\right)^{\frac{1}{2}} 
         + ||\bar{w}||\left( c_l\sum_{i=1}^{N_1}||\bar{\phi}_i||^2\right)^{\frac{1}{2}} 
         \hspace{1.2in} &&\\
  = -\nu \sum_{i=1}^{N_1}||\bar{\phi}_i(t)||^2+ \left( 2\sqrt{2N_1}c(\alpha)+1\right)
        ||\bar{w}||
        \left( c_l\sum_{i=1}^{N_1}||\bar{\phi}_i||^2\right)^{1/2} \hspace{2.2in} &&\\
  \leq -\frac{\nu}{2} \sum_{i=1}^{N_1}||\bar{\phi}_i(t)||^2
    +  \frac{\left( 2\sqrt{2N_1}c(\alpha)+1\right)^2c_l}{2\nu}
        ||\bar{w}||^2 \ \ \text{by Young's inequality}\hspace{1.6in} &&\\
  \leq -c_0\nu \lambda_1\frac{N_1 (N_1+1)}{4}
     + \frac{\left( 2\sqrt{2N_1}c(\alpha)+1\right)^2c_l}{2\nu}
        ||\bar{w}||^2 \hspace{2.9in}
\end{eqnarray*}
Hence,
{\small
\begin{equation*}
 \frac{1}{t}\int_0^t Tr (L_1(s)P_N(\bar{\Phi}_1(s),...,\bar{\Phi}_N(s))\,ds 
   \leq -c_0\nu \lambda_1\frac{N_1 (N_1+1)}{4} + 
      \frac{\left( 2\sqrt{2N_1}c(\alpha)+1\right)^2c_l}{2\nu}\frac{1}{t}
      \int_0^t||\bar{w}(s)||^2\,ds 
\end{equation*}
}

The estimate of the remaining trace follows \cite[p 60]{Markus}:
\[ \frac{1}{t}\int_0^t Tr(AP_N(\Phi_1^\bot (s),..., \Phi_N^\bot (s))\,ds \geq \frac{3}{10}c_0\lambda_1\nu N_2^{5/3}.
\]
Therefore,
\begin{eqnarray*}
   \frac{1}{t}\int_0^t Tr (L_c(s)P_N(\vec{\Phi}_1(s),...,\vec{\Phi}_N(s))\,ds 
    &\leq& -\nu \lambda_1 c_0\frac{N_1^2+N_2^{5/3}}{4} \nonumber \\
      & & \ \ \ \   +  \frac{\left( 2\sqrt{2N_1}c(\alpha)+1\right)^2c_l}{2\nu}\frac{1}{t}
      \int_0^t||\bar{w}(s)||^2\,ds.
\end{eqnarray*}

%-------------------
 \subsection{Dimensions}
%--------------------

Let $q_N(t)=  \frac{1}{t}\int_0^t Tr (L_c(s)P_N(\vec{\Phi}_1(s),...,\vec{\Phi}_N(s))\,ds$. Then,
\begin{eqnarray}
  q_N(t) &\leq& \nu \lambda_1 \left( -\frac{c_0}{4} (N_1^2+N_2^{5/3})
        +  \frac{\left(  2\sqrt{2N_1}c(\alpha)+1\right)^2c_l}{2\nu^2}\frac{1}{t\lambda_1}
      \int_0^t||\bar{w}(s)||^2\,ds\right) \nonumber \\
 \limsup_{t\rightarrow\infty} q_N(t) &\leq& \nu \lambda_1 \left( -\frac{c_0}{4} (N_1^2+N_2^{5/3})
   + \frac{\left(  2\sqrt{2N_1}c(\alpha)+1\right)^2c_l}{2\nu^2}\epsilon\right) , \nonumber
\end{eqnarray}
where $\epsilon = \nu\lambda_1 \limsup_{t\rightarrow \infty}\sup_{w_0\in X}\frac{1}{t\lambda_1}\int_0^t
||\bar{w}(s)||^2\,ds$
 with $X=\mathcal{A}_\alpha$.
 \smallskip 

To estimate $q_N$ in terms of $N$ ($=N_1+N_2$), a technical lemma is needed.

\begin{lemma}\label{lemma5.1}(\cite{Markus})
 Let $q\geq 0$. Then there exists a constant $c = c(q) > 0$ such that
 \[ x^q  + y^q \geq c (x+y)^q \]
 for all $x, y\geq 0$. 
\end{lemma}

{\bf Proof.} The result is true when $q=0$ with $c=1$ for all $x, y\geq 0$. It is also valid when $x=y=0$ for 
any $q>0$. Thus it will be sufficient to prove that, for any $q>0$,
there exists a constant $c>0$ such that 
\[ \frac{x^q+y^q}{(x+y)^q} = \frac{1+(y/x)^q}{(1+(y/x))^q} \geq c, \]
whenever $x, y>0$. By setting $z = y/x$, we can define a strictly positive and continuous function
on $(0, \infty )$
\[ f(z) = \frac{1+z^q}{(1+z)^q} \]
Notice that
\[ \lim_{z\rightarrow 0^+} f(z) = \lim_{z\rightarrow\infty}f(z) = 1, \]
and we can continuously extend $f$ on $(0, \infty)$ to $f_e$ on $[0,\infty )$.
Clearly $f_e(0) = 1$ and since $\lim_{z\rightarrow 0^+}f_e(z) = 1$,
there exists $z_0 > 0$ such that
\[ f_e(z) \geq \frac{1}{2} \]
whenever $z\geq z_0$. Since $[0, z_0]$ is compact, $f_e$ has an absolute minimum $m_0 = f_e(z_1)>0$
for some $z_1\in [0, z_0]$.
Choosing $c = c(q) = \min\{ 1/2, m_0\}$ proves the lemma. \ $\blacksquare$

We can make the constant $c$ in the Lemma \ref{lemma5.1} more precisely.

\begin{corollary}\label{cor5.2}
  \[ c(q) = \left\{ \begin{array}{ccc}
                     1 &\mbox{if} & 0< q \leq 1 \\
                     \frac{1}{2^{q-1}} &\mbox{if} & q>1.
                    \end{array} \right. 
  \]
  \end{corollary}
{\bf Proof.}  We can assume $z_0>1$ without loss of generality.
  To find a minimum $m_0>0$ of $f_e$ on $[0, z_0]$,
  \begin{eqnarray*}
   f_e^\prime(z) = 0 &\Leftrightarrow & qz^{q-1}(1+z)^q - q(1+z^q) (1+z)^{q-1} = 0 \\
     &\Leftrightarrow & q(1+z)^{q-1}(z^{q-1}-1) = 0 \\
     &\Leftrightarrow & z^{q-1}-1 = 0, \ \mbox{since\ $q>0$\ and}\ (1+z)^{q-1}\neq 0 \\
     &\Leftrightarrow & z^{q-1} = 1.
  \end{eqnarray*}
  
  \begin{itemize}
   \item \underline{Case 1. When $q=1$} \\
         $z^0=1$ for all $z\in (0, z_0)$; i.e., $f_e^\prime (z) = 0$ so that 
         $f_e(z)=1$ for all $z\in (0, z_0)$. Since $f_e(0)=1$ and $m_0=f_e(z_0)=1$
         by the continuity of $f_e$ on $[0, z_0]$. It implies $c=1$.
    \item \underline{Case 2. When $0<q<1$} \\
      \[ z^{q-1}=1 \Leftrightarrow \frac{1}{z^{1-q}} =1  \Leftrightarrow z^{1-q} = 1\]
      Taking the logarithm of both sides yields
      \[ (1-q)\ln z = 0 \Leftrightarrow z=1\ \ \mbox{since}\ 1-q>0. \]
      By the first derivative test in calculus, $f_e$ has a local maximum at $z=1$ on $[0, z_0]$, with
      $f_e(1) = 2/2^q>1$, so that the minimum must be $f_e(0)=1$, which is less than $f_e(z_0)$.
      Thus $c=1$.
   \item \underline{Case 3. When $q>1$} \\
    Similar to Case 2 with $q>1$. $f_e$ has an absolute minimum at $z=1$ on $[0, z_0]$, with
      $f_e(1) = 2/2^q = 2^{1-q}<1$.
      Thus, $c=2^{1-q}$.
  \end{itemize}
 This proves the corollary. 
  \ $\blacksquare$

By Corollary \ref{cor5.2},
\[ N_1^2+N_2^{5/3} \geq N_1^{5/3}+N_2^{5/3} \geq c\left(\frac{5}{3}\right) N^{5/3} = \left(\frac{1}{4}\right)^{1/3}N^{5/3}. \]
So,
\begin{eqnarray*}
 \lim_{t\rightarrow \infty}q_N(t) &\leq & \nu \lambda_1 \left( -d N^{5/3} +
           \frac{(2\sqrt{2}\sqrt{N} c(\alpha )+1)^2c_l}{2\nu^2}\epsilon\right) , \ \ \mbox{where}\ d=\frac{c_0}{4\sqrt[3]{c}}\\
        &=& \nu \lambda_1 \left( -d N^{5/3} +
           \frac{(8N c^2(\alpha )+4\sqrt{2}\sqrt{N}c(\alpha)+1)c_l}{2\nu^2}\epsilon\right) \\
           &\leq & \nu \lambda_1 \left( -d N^{5/3} +
           \frac{8\sqrt{2}N c(\alpha )+4\sqrt{2}Nc(\alpha)+1}{\nu^2}c_l\epsilon\right) \\
           & &  \mbox{since}\ c^2(\alpha) < 2 \ \mbox{and}\ N\geq 1\\
           &\leq & \nu \lambda_1 \left( -d N^{5/3} +
           \frac{24Nc(\alpha )+1}{\nu^2}c_l\epsilon\right)\\
         &\leq & \nu \lambda_1 \left( -d N^{5/3} +
           \frac{N(24c(\alpha )+1)}{\nu^2}c_l\epsilon\right)
\end{eqnarray*}

We want to find the smallest $N>0$ such that
\[
  -d N^{5/3} +
           \frac{N(24c(\alpha )+1)}{\nu^2}c_l\epsilon < 0.
\]
Setting the nonlinear equation
\begin{equation}\label{eq5.3}
  N^{5/3} = \frac{N(24c(\alpha )+1)}{d\nu^2}c_l\epsilon 
\end{equation}
yields 
\begin{eqnarray*}
  N &=& \frac{1}{\nu^3} \sqrt{ \frac{c_l^3\epsilon^3 (24c(\alpha)+1)^3}{d^3} } \\
  &=& \left(\frac{c_l}{d}\right)^{3/2} 
         \left( \frac{(24c(\alpha)+1)\epsilon}{\nu^2}\right)^{3/2} \\
  &=& \left(\frac{c_l}{d}\right)^{3/2} 
         (24c(\alpha)+1)^{3/2}\left(\frac{\epsilon}{\nu^2}\right)^{3/2}
\end{eqnarray*}
so that
\[
   N \geq \left(\frac{c_l}{d}\right)^{3/2} 
         (24c(\alpha)+1)^{3/2}\left(\frac{\epsilon}{\nu^2}\right)^{3/2}.
\]
Hence, 
\begin{eqnarray*}
 d_H(\mathcal{A}_\alpha) &<& \left(\frac{c_l}{d}\right)^{3/2} 
         (24c(\alpha)+1)^{3/2}\left(\frac{\epsilon}{\nu^2}\right)^{3/2} \\
         &<& \left(\frac{c_l}{d}\right)^{3/2} 
         (24c(\alpha)+1)^{3/2} \tilde{c}^{3/2}\left(\frac{\rho_V}{\nu^2}\right)^{3}\ \ \ \mbox{since}\ \ 
            \epsilon < \tilde{c}\rho^2_V \\
          &=& K(\alpha)\left(\frac{\rho_V}{\nu^2}\right)^{3}
\end{eqnarray*}
where 
\[ K(\alpha) =  \left(\frac{c_l}{d}\right)^{3/2} 
         (24c(\alpha)+1)^{3/2} \tilde{c}^{3/2} \rightarrow 
         \left(\frac{(24\sqrt{2}+1)c_l\tilde{c}}{d}\right)^{3/2} \equiv K_0\ \ \mbox{as}\ \alpha\rightarrow 0^+
\]

In particular, when $\alpha = 0$, the exact rotating Navier-Stokes equations
don't have the second term of the first inequality on page 17, and we get better estimate:
    \begin{equation*}
      d_H(\mathcal{A}_0) < \tilde{K}\left(\frac{\rho_{v_0}}{\nu_0}\right)^{6/5}
     \end{equation*}
This completes the proof of our main result Theorem \ref{main}.

%%%%%%%%%%%%%%%%%%%%%%%%%%%%%%%%%%%%%%%%%%%%%%%%%%%%%%%%%%%%%%%%

%%%%%%%%%%%%%%%%%%%%%%%%%%%%%%%%%%%%%%%%%%%%%%%%%%%%%%%%%%%%%%%%%%%

\end{document}